\documentclass[11pt]{amsart}

\usepackage{amssymb}
\usepackage{amsmath}
\usepackage{mathrsfs}
\usepackage{url}
\usepackage{verbatim}
\usepackage{multicol}
\input xy
\xyoption{all}
\usepackage{multirow}
\usepackage{mathdots}
\usepackage{hyperref}

\usepackage{graphicx}  
\usepackage{tikz}
\usetikzlibrary{decorations.pathreplacing}
\usepackage{tkz-graph}
\usetikzlibrary{arrows}
\pgfarrowsdeclare{mytip}{mytip}{%
  \setlength{\arrowsize}{1pt}  
  \addtolength{\arrowsize}{.5\pgflinewidth}  
  \pgfarrowsrightextend{0}  
  \pgfarrowsleftextend{-5\arrowsize}  
}{%
  \setlength{\arrowsize}{1pt}  
  \addtolength{\arrowsize}{.5\pgflinewidth}  
  \pgfpathmoveto{\pgfpoint{-5\arrowsize}{1.5\arrowsize}}  
  \pgfpathlineto{\pgfpointorigin}  
  \pgfpathlineto{\pgfpoint{-5\arrowsize}{-1.5\arrowsize}}  
  \pgfusepathqstroke  
}

\newtheorem{thm}{Theorem}
\newtheorem{lemma}[thm]{Lemma}

\newtheorem{cor}[thm]{Corollary}

\newtheorem{conj}[thm]{Conjecture}

\newtheorem{fact}[thm]{Fact}

\newtheorem{Definition}[thm]{Definition}

\newtheorem{Example}{Example}
\newenvironment{example}
  {\begin{Example}\rm}{\end{Example}}

\newcommand{\arxiv}[1]{\href{http://arxiv.org/abs/#1}{\texttt{arXiv:#1}}}

\title{Pattern Avoidance in Poset Permutations}

\author{Sam Hopkins}
\email{shopkins@mit.edu}
\address{Massachusetts Institute of Technology, Cambridge, MA, 02139}

\author{Morgan Weiler}
\email{morgan.w@berkeley.edu }
\address{University of California, Berkeley, Berkeley, CA,  94704}

\begin{document}

\begin{abstract}
  We extend the concept of pattern avoidance in permutations on a totally ordered set to pattern avoidance in permutations on partially ordered sets. The number of permutations on $P$ that avoid the pattern $\pi$ is denoted $Av_P(\pi)$. We extend a proof of Simion and Schmidt to show that $Av_P(132) \leq Av_P(123)$ for any poset $P$, and we exactly classify the posets for which equality holds.
\end{abstract}

\maketitle

\section{Motivation}

An inversion of a permutation $\sigma \in S_n$ is a pair of entries $a < b \in \{1,\ldots,n\}$ such that $b$ is to the left of $a$ in $\sigma$. As B\'{o}na \cite{bona} explains, classical pattern containment is ``a far-fetching generalization of [inversions of permutations] from pairs of entries to $k$-tuples of entries.'' A permutation $\sigma \in S_n$ is said to \emph{contain} a pattern $\pi \in S_k$ if some subsequence of $\sigma$ of length $k$ is order-isomorphic to $\pi$. Otherwise, we say $\sigma$ \emph{avoids} $\pi$. Thus an inversion in $\sigma$ is just a $21$ pattern that it contains.

We propose a similar generalization to permutations on posets. A total ordering of the elements of some poset $P$ that respects its partial order is called a \emph{linear extension} of $P$. It is possible also to think of a linear extension of $P$ as a permutation of the elements of $P$ that has no inversions. While the only classical permutation in $S_n$ that has no inversions is $12\ldots n$, in general there may be many linear extensions of $P$, and counting the number of such extensions is a difficult problem. Since linear extensions are a central object of study in order theory, it is natural to look at avoidance of more complicated patterns than $21$ in poset permutations. The field of pattern avoidance has blossomed in the past two decades, and in particular has expanded to include patterns in structures other than $S_n$ such as words~\cite{alon}~\cite{branden}, compositions and multiset permutations~\cite{atkinson}~\cite{albert}~\cite{savage}~\cite{heubach}, set partitions~\cite{klazar}~\cite{sagan}, ordered set partitions~\cite{godbole}, matchings~\cite{bloom}, et cetera. In~\cite{kitaev}, Kitaev studies classical permutation avoidance of patterns with incomparable elements. However, we do not believe pattern avoidance in poset permutations has been considered in the way we define below.

A central theme in the study of pattern avoidance is demonstrating relationships between the number of permutations that avoid different patterns of the same length. In this paper, we define a notion of poset pattern avoidance and demonstrate one non-trivial inequality between avoidance of two length-three patterns. In particular, our work builds on a dichotomy between `small-large-medium' and `small-medium-large' pattern avoidance in ordered structures. The following heuristic argument gives a reason for us to expect fewer `small-large-medium'-avoiders than `small-medium-large'-avoiders:
\begin{quote}
A random size three substructure of a given ordered structure is equally likely to be isomorphic to any of the size three order patterns. So among all structures of a given size, each size three pattern is contained an equal number of times in total. `Small-medium-large' substructures overlap one another better than do `small-large-medium' ones. So structures that contain `small-medium-large' substructures tend to contain a lot of them, as compared to `small-large-medium'-containing substructures. Thus, there should be more structures that avoid `small-medium-large' substructures all together.
\end{quote}
Our aim is to make this heuristic precise and understand when we have strict inequality versus equality. Indeed, our main result, Theorem~\ref{thm:123}, is an instance of this `small-large-medium' versus `small-medium-large' phenomenon. In~\S\ref{sec:apps} we give some applications of this theorem: specifically, we apply it to pattern avoidance in multiset permutations and in words (which are in fact a special case of poset permutation pattern avoidance), as well as to the avoidance of a new kind of pattern in classical permutations.

\section{Definitions and Basic Facts}

Let $P$ be a partially ordered set on $n$ elements. A permutation on $P$ is a bijection $\sigma \colon \{1,\ldots,n\} \to P$. We define $\sigma_i := \sigma(i)$ and think of $\sigma$ as an ordered list of the elements of $P$.  We will use the notation $\sigma = (\sigma_1,\ldots,\sigma_n)$ and we call $\sigma_i$ the \emph{entry} of $\sigma$ at position $i$. We will use $S_P$ to denote the set of permutations on $P$.

Let $Q$ be a poset on $k$ elements. A permutation $\sigma \in S_P$ is said to \emph{contain} a permutation $\pi \in S_Q$ (which we call a \emph{pattern}) if $\sigma$ has $k$ entries~$\sigma_{i_1}, \sigma_{i_2}, \ldots, \sigma_{i_k}$ so that $i_1 < i_2 < \cdots < i_k$, and for any $1 \leq a < b \leq k$, the order relation ($<$, $>$, or $\sim$ := ``incomparable'') between $\sigma_{i_a}$ and $\sigma_{i_b}$ is the same as the order relation between $\pi_a$ and $\pi_b$. Otherwise, we say $\sigma$ \emph{avoids} $\pi$. We call patterns $(\sigma_{i_1},\ldots,\sigma_{i_k})$ and $(\pi_1,\ldots,\pi_k)$ with the same order relations \emph{pattern isomorphic}. When considering permutations as patterns we will suppress the parentheses and write $\pi = \pi_1\pi_2\ldots\pi_n$. The most convenient notation to allow for arbitrary poset permutation patterns is to have the $\pi_i$ be sets with the partial order of containment. We will often use example permutations on $B_n$, the Boolean lattice on $n$ elements.

\begin{example}
Let $\sigma = (\{2,3\},\{2\},\{1,3\},\{1,2,3\},\{1\},\emptyset,\{1,2\},\{3\}) \in S_{B_3}$. Then $\sigma$ avoids $\emptyset\{1\}\{1,2\}$. However, $\sigma$ contains the pattern~$\{1\}\{3\}\{1,2\}$, as evidenced by the subsequence $(\{2\},\{1,3\},\{1,2\})$. $\square$
\end{example}

We will also consider patterns on chains within partially ordered sets in terms of their representation as patterns from the canonically totally ordered set $[k] := \{1,2,\ldots,k\}$. As usual, $S_k$ denotes the set of permutations on $[k]$.

\begin{example}
The $\sigma$ from the previous example avoids $123$, which is pattern isomorphic to $\emptyset\{1\}\{1,2\}$. $\square$
\end{example}

We will use both forms of notation but it will always be clear which notation is being used because of the set brackets. Note that $123 \neq \{1\}\{2\}\{3\}$.

For a permutation $\sigma$ on any poset $P$, define the \emph{reverse} of $\sigma$ to be the permutation~$(\sigma_n,\sigma_{n-1},\ldots,\sigma_1)$. For a permutation $\sigma$ on $P$, define the \emph{dual} of $\sigma$ to be the same list of entries $(\sigma_1,\ldots,\sigma_n)$ considered as a permutation on the dual poset of $P$.  We state the following without proof.

\begin{fact}
If $\sigma \in S_P$ avoids $\pi \in S_Q$, the reverse of $\sigma$ avoids the reverse of $\pi$.
\end{fact}

\begin{fact}
If $\sigma \in S_P$ avoids $\pi \in S_Q$, the dual of $\sigma$ avoids the dual of~$\pi$.
\end{fact} 

Denote the number of permutations on $P$ that avoid $\pi$ by $Av_P(\pi)$. In general, we cannot expect to find exact values of $Av_P(\pi)$. The simplest non-trivial case is where $\pi$ is $12$ or $21$. A \emph{linear extension} of a poset is a total ordering of its elements consistent with its partial ordering. Plainly, $Av_P(12) = Av_P(21)$ counts the number of linear extensions of $P$. Counting such extensions is often very hard. For example, for the classic example of the Boolean lattice $B_n$ we have only asymptotic bounds on the number of linear extensions due to Brightwell and Tetali~\cite{brightwell}. However, we may hope to find relations between the $Av_P(\pi)$. If $Av_P( \pi) = Av_P( \pi')$, we say that the patterns $\pi$ and~$\pi'$ are \emph{Wilf equivalent} for $P$. For instance, the above facts establish that $\pi$ and the reverse of $\pi$ are Wilf equivalent for all $P$, and that~$\pi$ and the dual of $\pi$ are Wilf equivalent for self-dual posets $P$.

The only other length two poset permutation pattern is $\{1\}\{2\}$. However, avoidance of $\{1\}\{2\}$ is trivial: there are no $\{1\}\{2\}$-avoiding permutations on~$P$ unless $P$ is a chain, in which case every permutation in $S_P$ avoids~$\{1\}\{2\}$.

\section{132- versus 123-avoiding poset permutations: An Injection}

A foundational result of classical permutation pattern avoidance is that there are as many $123$-avoiding permutations in $S_n$ as there are $132$-avoiding permutations, which when combined with the trivial observations about reverses and duals (or `complements') shows that $\sigma$ and $\sigma'$ are classically Wilf equivalent for all $\sigma,\sigma' \in S_3$. This has been proved in many ways, for instance in \cite[pp.~242-243]{knuth1} and \cite[pp.~60-61]{knuth2}, \cite{richards}, \cite{krattenthaler}, et cetera, but perhaps the most elegant proof is due to Simion and Schmidt~\cite{simion}. We extend the proof of Simion and Schmidt and exactly classify those posets for which $123$ and $132$ are Wilf equivalent patterns.

%edit: from here to the theorem statement was contained in the proof. I moved it above the proof re the referee's suggestion but you may want to do some more tinkering.

An important tool from classical pattern avoidance is the use of left-to-right minima to demonstrate bijections. An entry of $\sigma \in S_n$ is called a \emph{left-to-right minimum} (LRM) if it is less than every entry to its left. Fixing the positions of the LRM of $\sigma \in S_n$, there is exactly one way to fill in the remaining entries to yield a permutation with those LRM which avoids $132$. There is also exactly one way to fill in the entries to yield a permutation with those LRM which avoids $123$. Proving this rigorously gives Simion and Schmidt's bijection between $132$- and $123$-avoiding permutations.

\begin{example}
With $\rho = 67341258 \in S_8$, the LRM of $\rho$ are $6$, $3$, and $1$ in positions $1$, $3$ and $5$.  The only permutation in $S_8$ with those LRM in those positions that avoids $132$ is $\rho$, and the only permutation in $S_8$ with those LRM in those positions that avoids $123$ is $68371542$. $\square$
\end{example}

\begin{thm}
We have $Av_P(132) \leq Av_P(123)$ for any poset $P$, with strict inequality if and only if $P$ contains one of $Q_1$, $Q_2$, or $Q_3$ below as an induced subposet:
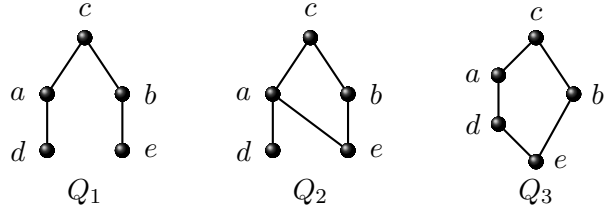
\begin{center}
\begin{tikzpicture}
\SetVertexMath
\GraphInit[vstyle=Art]
\SetUpVertex[MinSize=3pt]
\SetVertexLabel
\tikzset{VertexStyle/.style = {shape = circle,shading = ball,ball color = black,inner sep = 2pt}}
\SetUpEdge[color=black]

\Vertex[LabelOut,Lpos=180,Ldist=.05cm,x=0,y=0.25]{d}
\Vertex[LabelOut,Lpos=0,Ldist=.05cm,x=1,y=0.25]{e}
\Vertex[LabelOut,Lpos=90,Ldist=.05cm,x=0.5,y=1.75]{c}
\Vertex[LabelOut,Lpos=180,Ldist=.05cm,x=0,y=1]{a}
\Vertex[LabelOut,Lpos=0,Ldist=.05cm,x=1,y=1]{b}
\Edges(e,b,c,a,d)

\draw (0.5,-0.3) node{$Q_1$};

\Vertex[LabelOut,Lpos=180,Ldist=.05cm,x=3,y=0.25]{d}
\Vertex[LabelOut,Lpos=0,Ldist=.05cm,x=4,y=0.25]{e}
\Vertex[LabelOut,Lpos=90,Ldist=.05cm,x=3.5,y=1.75]{c}
\Vertex[LabelOut,Lpos=180,Ldist=.05cm,x=3,y=1]{a}
\Vertex[LabelOut,Lpos=0,Ldist=.05cm,x=4,y=1]{b}
\Edges(a,e,b,c,a,d)

\draw (3.5,-0.3) node{$Q_2$};

\Vertex[LabelOut,Lpos=180,x=6,y=0.6]{d}
\Vertex[LabelOut,Lpos=0,x=6.5,y=0.1]{e}
\Vertex[LabelOut,Lpos=90,x=6.5,y=1.75]{c}
\Vertex[LabelOut,Lpos=180,x=6,y=1.25]{a}
\Vertex[LabelOut,Lpos=0,x=7,y=1]{b}
\Edges(d,e,b,c,a,d)

\draw (6.5,-0.3) node{$Q_3$};

\end{tikzpicture}
\end{center}
\label{thm:123}
\end{thm}

%warn them about long proof

\noindent \emph{Proof}. We extend the concept of left-to-right minimum for permutations on a poset $P$ with what we call a \emph{left-to-right minimal element} (LRME): an LRME of $\sigma \in S_P$ is an entry $\sigma_i$ such that there is no $j < i$ with $\sigma_j < \sigma_i$. Informally, an LRME is less than or incomparable to every entry preceding it. Unlike in the classical permutations case, when fixing the positions of the LRME of a poset permutation there may be more than one way to fill in the remaining entries to yield a permutation which avoids either $132$ or~$123$.

\begin{example}
If $\sigma = (\{2,3\},\{1\},\{1,2\},\{2\},\emptyset,\{3\},\{1,3\},\{1,2,3\}) \in S_{B_3}$, the LRME of $\sigma$ are $\{2,3\}$, $\{1\}$, $\{2\}$, and $\emptyset$ in positions $1$, $2$, $4$, and $5$. Note that $\sigma$ avoids $132$ but the following other elements of $S_{B_3}$ have the same LRME in the same positions as $\sigma$ and also avoid $132$:
\begin{itemize}
	\item $(\{2,3\},\{1\},\{1,3\},\{2\},\emptyset,\{1,2\},\{3\},$ $\{1,2,3\})$;
	\item $(\{2,3\},\{1\},\{1,3\},\{2\},\emptyset,\{3\},\{1,2\},\{1,2,3\})$.
\end{itemize}
Similarly, the following elements of $S_{B_3}$ all have the same LRME in the same positions and avoid $123$:
\begin{itemize}
	\item $(\{2,3\},\{1\},$ $\{1,2,3\},\{2\},\emptyset,\{1,2\},\{1,3\},\{3\})$; 
	\item $(\{2,3\},\{1\},\{1,2,3\},\{2\},\emptyset, \{1,3\}, \{1,2\}, \{3\})$;
	\item $(\{2,3\},\{1\},\{1,2,3\},\{2\},\emptyset,\{1,3\},\{3\},\{1,2\})$. $\square$
\end{itemize}
\end{example}

By an \emph{LRME set} we mean a list of elements $x_1, \ldots, x_k$ from $P$ along with a list of corresponding positions $1 \leq \mu_1 < \cdots < \mu_k \leq n$. Call an LRME set $X$ \emph{admissible} if there is some permutation $\sigma$ whose LRME are exactly~$x_1, \ldots, x_k$ in positions $\mu_1, \ldots, \mu_k$, and in this case we say $X$ is the LRME set of $\sigma$. Fix some admissible LRME set $X$. How many $\sigma$ have $X$ as their LRME set and avoid $132$?

Denote the positions not among the $\mu_i$ by \mbox{$1 \leq \nu_1 < \cdots < \nu_l \leq n$}. Let~$P'$ be the induced subposet of $P$ on the elements that are not among the~$x_i$. First note that since $\sigma_{\nu_k}$ is not an LRME of $\sigma$ for any $1\leq k\leq l$, $\sigma_{\nu_k}$ is greater than some $\sigma_{\mu_j}$, where $\mu_j < \nu_k$. So, for $y \in P'$, define $\omega(y)$ to be the the smallest $i$ such that there exists $\sigma_{\mu_j} < y$ for some $\mu_j < \nu_i$. Call a permutation $\sigma'$ on $P'$ \emph{$\omega$-legal} when it obeys the condition that $\omega(\sigma'_i) \leq i$ for all $i$, and let $\Lambda^{\omega}$ be the set of $\omega$-legal permutations in $S_{P'}$. Then, $\sigma$ has $X$ as its LRME set exactly when the subsequence $(\sigma_{\nu_1}, \ldots, \sigma_{\nu_l})$ is in $\Lambda^{\omega}$.

\begin{example}
With $\sigma = (\{2,3\},\{1\},\{1,2\},\{2\},\emptyset,\{3\},\{1,3\},\{1,2,3\}) \in S_{B_3}$ as in the previous example, $\mu_1 = 1$, $\mu_2 = 2$, $\mu_3 = 4$, $\mu_4 = 5$, $\nu_1 = 3$, $\nu_2 = 6$, $\nu_3 = 7$ and $\nu_4 = 8$. We have $\omega(\{1,3\}) = \omega(\{1,2\}) = \omega(\{1,2,3\}) = 1$ since each of these elements is greater than $\sigma_{\mu_2}$, while $\omega(\{3\}) = 2$. Note that~$\omega$ values refer to the $\nu$ indices, \emph{not} to the positions of the elements in~$\sigma$. $\square$
\end{example}

If $\sigma$ contains a $132$ pattern then it contains a $132$ pattern that consists of an LRME followed by two non-LRME, since if the element acting as the~1 is not an LRME then the rightmost LRME to the left of it that it is greater than will also make a 132 pattern with the same elements acting as~$3$ and~$2$. Suppose we fill each $\nu_i$ from left to right by choosing an $\omega$-legal, unchosen element of $P'$ to occupy this position. If we ever choose $z$ when $y$ is also an $\omega$-legal choice and $y < z$, we will contain a $132$ (with the $32$ being $zy$ and the $1$ being the LRME they both are greater than). If, on the other hand, we always choose a minimal $\omega$-legal element, we will avoid $132$; in any $132$ pattern $xzy$, with $x$ an LRME and $z$ and $y$ non-LRME, $y$ was an~$\omega$-legal choice for the position $z$ occupies. Let $\Lambda^{\omega}_{\mathrm{min}} \subset \Lambda^{\omega}$ be those~$\omega$-legal permutations $\sigma'$ for which as we fill in the entries from left to right we always choose a minimal element among the $\omega$-legal choices. Then $\sigma$ has~$X$ as its LRME set and avoids $132$ exactly when $(\sigma_{\nu_1}, \ldots, \sigma_{\nu_l})$ is in $\Lambda^{\omega}_{\mathrm{min}}$.

Similarly, $\sigma$ avoids $123$ if and only if as we fill in the $\nu_i$ from left to right we always choose a maximal element among $\omega$-legal choices.  Let $\Lambda^{\omega}_{\mathrm{max}} \subset \Lambda^{\omega}$ be those $\omega$-legal permutations $\sigma'$ for which as we fill in the entries from left to right we always choose a maximal element among the $\omega$-legal choices. Then $\sigma$ has $X$ as its LRME set and avoids $123$ exactly when $(\sigma_{\nu_1}, \ldots, \sigma_{\nu_l})$ is in $\Lambda^{\omega}_{\mathrm{max}}$. In order to complete the proof, we need to show $|\Lambda^{\omega}_{\mathrm{min}}| \leq |\Lambda^{\omega}_{\mathrm{max}}|$. The following lemma, embedded within the proof of Theorem~\ref{thm:123}, will give us just that. Specifically, fixing some $X$ and defining $\omega$ in terms of $X$, we will apply the following lemma to this $\omega$.

\begin{lemma}
Let $P$ be a poset on $n$ elements. Define a partial order $\preceq$ on~$S_P$ whereby for permutations $\sigma, \pi \in S_P$ we say $\pi \prec \sigma$ if there exists $1 \leq j \leq n$ such that $\pi_j < \sigma_j$ and $\sigma_i = \pi_i$ for all $i < j$.

Let $\omega \colon P  \to \{1,\ldots,n\}$ be a labeling function such that
\begin{enumerate}
\item the number of elements $x$ with $\omega(x) \leq i$ is greater than or equal to $i$ for all $i = 1, \ldots, n$, and, \label{cond:nonempty}
\item if $x > y$, then $\omega(x) \leq \omega(y).$ \label{cond:decreasing}
\end{enumerate}
Call a permutation $\sigma \in S_P$ \emph{$\omega$-legal} if $\omega(\sigma_i) \leq i$ for all $i = 1, \ldots, n$. Let~$\Lambda^{\omega}$ be the set of $\omega$-legal permutations on $P$ and define 
\[\Lambda^{\omega}_{\mathrm{max}} := \{\sigma \in \Lambda^{\omega}\colon \nexists \; \sigma' \in \Lambda^{\omega} \textrm{ such that } \sigma \prec \sigma'\}\] 
\[\Lambda^{\omega}_{\mathrm{min}} := \{\sigma \in \Lambda^{\omega}\colon \nexists \; \sigma' \in \Lambda^{\omega} \textrm{ such that } \sigma' \prec \sigma\}.\]
Then there is an injection
\[ \phi\colon \Lambda^{\omega}_{\mathrm{min}} \hookrightarrow \Lambda^{\omega}_{\mathrm{max}}.\] 
Further, $\phi$ is a bijection if there does not exist $x,y,z \in P$ with 
\begin{itemize}
	\item $x < z$; 
	\item $y < z$; 
	\item $x \sim y$;
	\item $\omega(x) < \omega(y)$. 
\end{itemize}
\label{lem:labeling}
\end{lemma}

\noindent \emph{Proof.} Condition (\ref{cond:nonempty}) on $\omega$ merely guarantees that $\Lambda^{\omega}$ is nonempty (and consequently $\Lambda^{\omega}_{\mathrm{min}}$ and $\Lambda^{\omega}_{\mathrm{max}}$ are also nonempty).

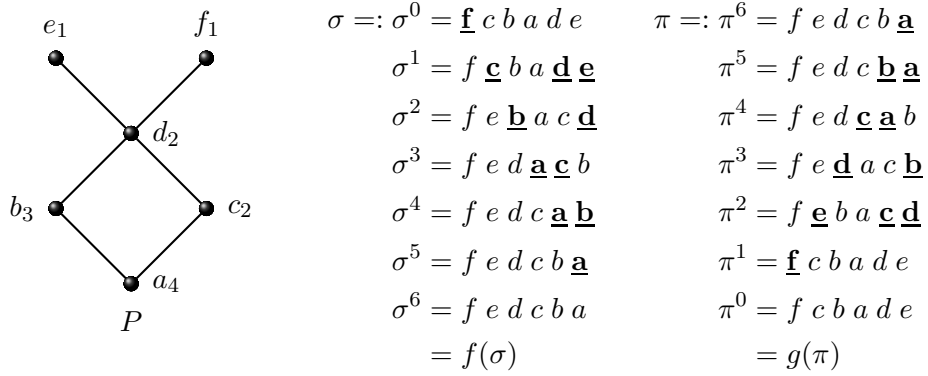
\begin{figure}[ht]
\begin{multicols}{3}
\begin{tikzpicture}
\SetVertexMath
\GraphInit[vstyle=Art]
\SetUpVertex[MinSize=3pt]
\SetVertexLabel
\tikzset{VertexStyle/.style = {shape = circle,shading = ball,ball color = black,inner sep = 2pt}}
\SetUpEdge[color=black]
\Vertex[LabelOut,Lpos=0,Ldist=.05cm,x=1,y=0]{a_4}
\Vertex[LabelOut,Lpos=180,Ldist=.05cm,x=0,y=1]{b_3}
\Vertex[LabelOut,Lpos=0,Ldist=.05cm,x=2,y=1]{c_2}
\Vertex[LabelOut,Lpos=0,Ldist=.05cm,x=1,y=2]{d_2}
\Vertex[LabelOut,Lpos=90,Ldist=.05cm,x=0,y=3]{e_1}
\Vertex[LabelOut,Lpos=90,Ldist=.05cm,x=2,y=3]{f_1}
\Edges(f_1,d_2,b_3,a_4,c_2,d_2,e_1)
\draw (1,-0.5) node{$P$};
\end{tikzpicture}
\medskip
\begin{align*}
\sigma =: \sigma^0 &= \underline{\bf f} \; c \; b \; a \; d \; e \\
\sigma^1 &= f \; \underline{\bf c} \; b \; a \; \underline{\bf d} \; \underline{\bf e} \\
\sigma^2 &= f \; e \; \underline{\bf b} \; a \; c \; \underline{\bf d} \\
\sigma^3 &= f \; e \; d \; \underline{\bf a} \; \underline{\bf c} \; b \\
\sigma^4 &= f \; e \; d \; c \; \underline{\bf a} \; \underline{\bf b} \\
\sigma^5 &= f \; e \; d \; c \; b \; \underline{\bf a} \\
\sigma^6 &= f \; e \; d \; c \; b \; a\\
&=f(\sigma)
\end{align*}
\medskip
\begin{align*}
\pi =: \pi^6 &= f \; e \; d \; c \; b \; \underline{\bf a} \\
\pi^5 &= f \; e \; d \; c \; \underline{\bf b} \; \underline{\bf a} \\
\pi^4 &= f \; e \; d \; \underline{\bf c} \; \underline{\bf a} \; b \\
\pi^3 &= f \; e \; \underline{\bf d} \; a \; c \; \underline{\bf b} \\
\pi^2 &= f \; \underline{\bf e} \; b \; a \; \underline{\bf c} \; \underline{\bf d} \\
\pi^1 &= \underline{\bf f} \; c \; b \; a \; d \; e \\
\pi^0 &= f \; c \; b \; a \; d \; e \\
&=g(\pi)
\end{align*}
\end{multicols}
\caption{Example run of $f(\sigma)$ algorithm and $g(\pi)$ algorithm. The subscript of each element $x$ in the Hasse diagram of $P$ is $\omega(x)$. For each $\sigma^{i}$ and $\pi^{i}$ in the algorithm, the entries at $\alpha_0, \ldots, \alpha_k$ and $\beta_0,\ldots,\beta_l$, respectively, are bold and underlined.}\label{fig:algorithmex}
\end{figure}

We now define a function $f\colon \Lambda^{\omega} \to \Lambda^{\omega}_{\mathrm{max}}$, whose restriction to $\Lambda^{\omega}_{\mathrm{min}}$ will be the $\phi$ we are looking for. The following algorithm defines $f$.  Let $\sigma \in \Lambda^{\omega}$. We will build a series of permutations $\sigma^{0}, \sigma^{1}, \ldots, \sigma^{n}$. Initialize $\sigma^{0} := \sigma$. When we are done, $f(\sigma)$ will be defined as $\sigma^{n}$. We recursively define $\sigma^{i+1}$ from $\sigma^{i}$:
\begin{enumerate} 
\item Mark position $i+1$.
\item Consider each position $j$ with $i+2 \leq j \leq n$ from left to right. If the entry $\sigma_j$ at the corresponding position is greater than the entry at the last marked position, mark $j$. 
\item Let $\alpha_0, \ldots, \alpha_k$ be the list of marked positions in the order they were marked. 
\item Set $\sigma^{i+1}_{\alpha_0} = \sigma^{i}_{\alpha_k}, \sigma^{i+1}_{\alpha_1} = \sigma^{i}_{\alpha_0}, \ldots, \sigma^{i+1}_{\alpha_k} = \sigma^{i}_{\alpha_{k-1}}$. For all other positions, set the entry of $\sigma^{i+1}$ to be the same as $\sigma^{i}$. In other words, let~$\sigma^{i+1} :=  \sigma^{i} \circ \gamma$, where $\gamma \in S_n$ is the cycle $(\alpha_k,\alpha_{k-1},\ldots,\alpha_1,\alpha_0)$.
\end{enumerate}
Figure~\ref{fig:algorithmex} gives an example of one run of the algorithm.

We claim that for any $\sigma \in \Lambda^{\omega}$, the permutation $f(\sigma)$ is $\omega$-legal and in particular $f(\sigma) \in \Lambda^{\omega}_{\mathrm{max}}$. Further we claim that there exists some function~$g\colon \Lambda^{\omega}_{\mathrm{max}} \to \Lambda^{\omega}$ such that $g(f(\sigma)) = \sigma$ for all $\sigma \in \Lambda^{\omega}_{\mathrm{min}}$, i.e., that $\phi$ is injective.

We move an element leftward at any step in the algorithm only if it was greater than the element previously occupying the position to which it moves. Thus if $\sigma^{i}$ was an $\omega$-legal permutation then so is $\sigma^{i+1}$ because condition (\ref{cond:decreasing}) on $\omega$ gives $\omega(\sigma^{i+1}_j) \leq \omega(\sigma^{i}_j) \leq j$ for all $j$. So $\omega$-legality is maintained at every step of the algorithm, and therefore $f(\sigma)$ is an $\omega$-legal permutation.

Clearly $f(\sigma) \in  \Lambda^{\omega}_{\mathrm{max}}$ because at every step in the algorithm we set $\sigma^{i+1}_{i}$ to be a maximal element among elements to the right of position $i$ in $\sigma^{i}$ and after this step the entry at position $i$ never changes.

To show that $\phi$ is injective, we construct its left inverse. The following algorithm defines $g\colon \Lambda^{\omega}_{\mathrm{max}} \to \Lambda^{\omega}$. Let $\pi \in \Lambda^{\omega}_{\mathrm{max}}$. We will build a series of permutations $\pi^{n}, \pi^{n-1}, \ldots, \pi^{0}$. Initialize $\pi^{n} := \pi$. When we are done, $g(\pi)$ will be defined as $\pi^{0}$. We recursively define $\pi^{i-1}$ from $\pi^{i}$:
\begin{enumerate}
\item Mark position $i$. 
\item Consider each position $j$ with $n\geq j \geq i+1$ from right to left. If the entry $\pi_j$ at the corresponding position is less than the entry at the last marked position, and if $\omega(\pi_j) \leq i$, mark $j$. 
\item Let $\beta_0, \beta_l, \beta_{l-1}, \ldots, \beta_1$ be the list of marked positions in the order they were marked. 
\item Set $\pi^{i-1}_{\beta_0} = \pi^{i}_{\beta_{1}}, \pi^{i-1}_{\beta_{1}} = \pi^{i}_{\beta_{2}}, \ldots, \pi^{i-1}_{\beta_l} = \pi^{i}_{\beta_{0}}$. For all other positions, set the entry of $\pi^{i-1}$ to be the same as $\pi^{i}$. That is, let ~$\pi^{i-1} := \pi^{i} \circ \gamma$, where $\gamma \in S_n$ is the cycle $(\beta_0,\beta_{1},\ldots,\beta_{l-1},\beta_l)$.
\end{enumerate}

We now show that $g(f(\sigma)) = \sigma$ for any $\sigma \in \Lambda^{\omega}_{\mathrm{min}}$. The proof that follows is technical but necessary. Call a permutation $\sigma$ on $P$ \emph{$i$-minimal} if for any~$k > j > i$ with $\sigma_k < \sigma_j$, we have $\omega(\sigma_k) > j$. This property will be useful for showing that the $(n-i)$-th step of the $g$ algorithm undoes the $i$th step of the $f$ algorithm because during these steps we consider only the entries in positions $i+1$ to $n$. First we claim that each $\sigma^{i}$ in the $f(\sigma)$ algorithm is $i$-minimal. We prove this by induction. The case $i = 0$ holds because $0$-minimality is equivalent to $\sigma$ belonging to $\Lambda^{\omega}_{\mathrm{min}}$. An element moves leftward at the $i$th step only if it moves into position $i$. There is no step where element $x$ moves rightward past an element $y$ greater than it, because such a $y$ would be a part of that step's cycle. Thus if $x$ is to the left of $y$ in $\sigma^{i}$ with $x < y$, then $x$ is to the right of $y$ in $\sigma^{i+1}$ only if $y$ moves into position~$i$ during this step. So the claim follows by induction. %maybe be a little more careful

Set $\pi := f(\sigma)$ and consider the $\pi^{i}$ from the $g(\pi)$ algorithm. We prove by downward induction that $\pi^{i} = \sigma^{i}$ for all $i$. The case $i = n$ holds by definition. Assume that $\pi^{i} = \sigma^{i}$. We will show $\pi^{i-1} = \sigma^{i-1}$. Let~$\alpha_0, \ldots, \alpha_k$ be as defined in the $f(\sigma)$ algorithm at the step where we go from $\sigma^{i-1}$ to~$\sigma^{i}$. Let~$\beta_0, \ldots, \beta_l$ be as defined in the $g(\pi)$ algorithm at the step where we go from $\pi^{i}$ to $\pi^{i-1}$. Of course, $\alpha_0 = \beta_0$. Further, we have that \mbox{$\alpha_k = \beta_l$}. To see this, suppose that in the $g(\pi)$ algorithm, as we consider $j$ with~\mbox{$n \geq j \geq i + 1$} from right to left, we mark a position $\beta_l$ before $\alpha_k$; that is, suppose \mbox{$\beta_l < \alpha_k$}. Then, $\pi^{i}_{\beta_l} < \pi^{i}_{\beta_0}$ and~$\omega(\pi^{i}_{\beta_l}) \leq i \leq \alpha_k$. But $\sigma^{i-1}_{\alpha_k} = \sigma^{i}_{\alpha_0} = \pi^{i}_{\beta_0}$ and also~\mbox{$\sigma^{i-1}_{\beta_l} = \sigma^{i}_{\beta_l} = \pi^{i}_{\beta_l}$}, so we have $i - 1 < \alpha_k < \beta_l$ such that~\mbox{$\sigma^{i-1}_{\alpha_k} < \sigma^{i-1}_{\beta_l}$} and~$\omega(\sigma^{i-1}_{\beta_l}) \leq \alpha_k$, a contradiction of the $(i-1)$-minimality of $\sigma^{i-1}$. It cannot be that~\mbox{$\beta_l < \alpha_k$}: we definitely mark $\alpha_k$ when we come to it because~$\pi^{i}_{\alpha_k} = \sigma^{i}_{\alpha_k} = \sigma^{i-1}_{\alpha_{k-1}} < \sigma^{i-1}_{\alpha_{k}} = \sigma^{i}_{\alpha_0} = \pi^{i}_{\beta_0}$ and $\pi^{i}_{\alpha_k} = \sigma^{i-1}_{\alpha_{k-1}} > \sigma^{i-1}_{\alpha_{0}}$, which means $\omega(\pi^{i}_{\alpha_k}) \leq \omega(\sigma^{i-1}_{\alpha_{0}}) \leq i$. So $\alpha_k = \beta_l$. Applying this argument again gives $\alpha_{k-1} = \beta_{l-1}$, and so on; it also proves $k = l$. Thus we have that~$\alpha_i = \beta_i$ for all $i$. Then,
\begin{align*}
\pi^{i-1} &= \pi^{i} \circ (\beta_0,\beta_{1},\ldots,\beta_{k-1},\beta_k) \\
&= \sigma^{i} \circ (\alpha_0, \ldots, \alpha_k) \\
&= \sigma^{i-1}  \circ (\alpha_k,\ldots,\alpha_0)  \circ (\alpha_0, \ldots, \alpha_k) \\
&= \sigma^{i-1}.
\end{align*}
That $g(f(\sigma)) = \sigma$ follows by induction.

In fact, $g$ is also injective. That is, we have $f(g(\pi)) = \pi$ for all~\mbox{$\pi \in \Lambda^{\omega}_{\mathrm{max}}$}, which can be proved in a manner very similar to the above proof that we have~\mbox{$g(f(\sigma)) = \sigma$} for any $\sigma \in \Lambda^{\omega}_{\mathrm{min}}$. Thus $\phi$ is a bijection between $\Lambda^{\omega}_{\mathrm{min}}$ and~$\Lambda^{\omega}_{\mathrm{max}}$ if and only if $g(\pi) \in \Lambda^{\omega}_{\mathrm{min}}$ for every $\pi \in \Lambda^{\omega}_{\mathrm{max}}$. Suppose there exists $\pi \in \Lambda^{\omega}_{\mathrm{max}}$ such that $g(\pi)$ is not in $\Lambda^{\omega}_{\mathrm{min}}$. Let $\pi^{n}, \ldots, \pi^{0}$ be as defined in the $g(\pi)$ algorithm and let $i$ be the largest value such that $\pi^{i}$ is not~$i$-minimal. There must be such an $i$ because $\pi^{0} = g(\pi)$ is not $0$-minimal as it is not in $\Lambda^{\omega}_{\mathrm{min}}$. Also, $i$ must be less than $n$ because any permutation is trivially $n$-minimal. Then consider the step of the algorithm that takes us from $\pi^{i+1}$ to $\pi^{i}$. If, as we were marking positions from $n$ to $i+1$, we marked each position whose entry was less than the entry of the last marked position, we would maintain $i$-minimality. So it must be that we skip over some entry $y$ because $\omega(y) > i$. Let $z$ be the entry of the position we had marked before considering $y$ and let $x$ be the entry of the next position we mark after $y$. There must be some such $x$ so that $z$ moves to the left of $y$; in fact $x$ must be in position $j$ in $\pi^{i+1}$ with $\omega(y) \leq j$ so that $z$ moving into position $j$ violates $i$-minimality. Then $y < z$, $x < z$ and $\omega(x) < \omega(y)$. Of course $y$ is not greater than $x$, but further $x$ is not greater than $y$ as $x$ and~$y$ would then violate the $(i+1)$-minimality of $\pi^{i+1}$. So $x \sim y$ as claimed. $\square$

We now finish the proof of Theorem~{\ref{thm:123}}. The $\omega$ defined earlier in the proof of Theorem~{\ref{thm:123}} by choosing some LRME set~$X$ obeys conditions~(\ref{cond:nonempty}) and~(\ref{cond:decreasing}) from Lemma~{\ref{lem:labeling}} and the~$\Lambda^{\omega}_{\mathrm{min}}$ and~$\Lambda^{\omega}_{\mathrm{max}}$ defined earlier are the same as those in Lemma~{\ref{lem:labeling}}. Thus the injection $\phi$ from $\Lambda^{\omega}_{\mathrm{min}}$ to $\Lambda^{\omega}_{\mathrm{max}}$ gives rise to an injection from the set of permutations $\sigma \in S_P$ that have $X$ as their LRME set and avoid $132$ and those $\sigma$ that have $X$ as their LRME set and avoid $123$. By summing over all LRME sets $X$, we conclude that~$Av_P(132) \leq Av_P(123)$.

If $Av_P(132) < Av_P(123)$, then there has to be an admissible LRME set~$X$ such that $|\Lambda^{\omega}_{\mathrm{min}}| < |\Lambda^{\omega}_{\mathrm{max}}|$. In this case, the injection $\phi$ from Lemma~{\ref{lem:labeling}} must not be a bijection, and so there must be elements $a,b,c \in P$ with 
\begin{itemize}
	\item $a < c$;
	\item $b < c$; 
	\item $a \sim b$; 
	\item $\omega(a) < \omega(b)$. 
\end{itemize}
But $\omega(a) < \omega(b)$ only if the leftmost LRME that $a$ is greater than, call it~$d$, is to the left of the leftmost LRME that $b$ is greater than, call it $e$. Then the induced subposet of $P$ on $\{a,b,c,d,e\}$ matches one of $Q_1$, $Q_2$, or $Q_3$.

On the other hand, if $P$ contains as an induced subposet any of $Q_1$, $Q_2$, or~$Q_3$, there is some set of LRME such that there are strictly more $123$-avoiding permutations with these LRME than $132$-avoiding permutations. If there exists $c' \in P$ with $c' > c$ then the induced subposet on $\{a,b,c',d,e\}$ will be the same as on $\{a,b,c,d,e\}$ and so without loss of generality we may assume $c$ is maximal. Consider the permutation $\sigma = (\theta_1,d,c,\theta_2,e,\theta_3,a,b)$, where $(\theta_1,d,\theta_2,e,\theta_3)$ is a $12$-avoiding subsequence of $\sigma$ containing all the elements of $P \setminus \{a,b,c\}$ (here $\theta_1$, $\theta_2$, and $\theta_3$ are themselves permutations). Such a $12$-avoiding subsequence exists because $e$ is not greater than $d$. Consider all permutations with the same LRME set as $\sigma$. The non-LRME elements are $a$, $b$, and $c$, with $\omega(c) = \omega(a) = 1$ and $\omega(b) = 2$. It is easily seen that~$\Lambda^{\omega}_{\mathrm{min}} = \{abc\}$ while $\Lambda^{\omega}_{\mathrm{max}} = \{cab,cba\}$, so $Av_P(132) < Av_P(123)$. $\square$

\section{Applications of the 132 versus 123 result} \label{sec:apps}

\subsection{Multiset permutations and words} If we take $P = [n]$ then Theorem~\ref{thm:123} recaptures the result that there are the same number of permutations in $S_n$ that avoid $132$ as $123$. Of course, in this case the bijection is the same as that of Simion and Schmidt. Slightly more generally, we can recapture the analogous statement for pattern avoidance in multiset permutations or in words. The following is well-known (see~\cite[Theorem 4.7]{heubach}).

\begin{cor}
For any vector $\vec{a} = (a_1,\ldots,a_n)$, the number of permutations of the multiset $1^{a_1}2^{a_2}\ldots n^{a_n}$ that avoid $132$ is the same as that avoid $123$.
\end{cor}
\noindent {\emph Proof}: Consider the poset $P_{\vec{a}}$ whose elements are $i_j$ with $1 \leq i \leq n$ and~$1 \leq j \leq a_i$, where $i_j \leq i'_{j'}$ if and only if $i \leq i'$. Theorem~\ref{thm:123} tells us that $Av_P(132) = Av_P(123)$. Set $N := a_1!a_2!\cdots a_n!$. There is an $N$-to-one surjection from $S_{P_{\vec{a}}}$ to the set of permutations of the multiset $1^{a_1}2^{a_2}\ldots n^{a_n}$ that simply forgets subscripts. This map clearly preserves containment of the patterns $132$ and $123$. $\square$

By considering all vectors $(a_1,\ldots,a_n)$ with $a_1+\cdots+a_n = \ell$ for some fixed~$\ell$, we can obtain from the previous corollary an analogous statement about the number of words of length $\ell$ in the alphabet $\{1,\ldots,n\}$ that avoid~$132$ versus $123$. That result is also well-known (see~\cite[Theorem 4.8]{branden}).

\subsection{Gap patterns in classical permutations} Another application of Theorem~\ref{thm:123} is to classical permutations, but concerning the avoidance of a new, nonstandard kind of pattern. For $\ell \geq 0$, let us say that~$\sigma \in S_n$ contains the \emph{gap pattern} $1-_\ell2-_\ell3$ (respectively, $1-_\ell3-_\ell2$) if there exist three indices $1\leq i \leq j\leq k \leq n$ with $i + \ell < j$ and $j + \ell < k$ such that~$\sigma_i < \sigma_j < \sigma_k$ (resp.,~$\sigma_i < \sigma_k < \sigma_j$). In other words, we require that the representatives of the pattern be sufficiently far apart from one another. Denote by $Av_n(1-_\ell2-_\ell3)$ (respectively, $Av_n(1-_\ell3-_\ell2)$) the number of permutations in $S_n$ that avoid the gap pattern $1-_\ell2-_\ell3$ (resp., $1-_\ell3-_\ell2$).

\begin{cor}
For $\ell \geq 0$ and $n \geq 1$ we have $Av_n(1-_\ell3-_\ell2) \leq Av_n(1-_\ell2-_\ell3)$ with strict inequality if and only if $\ell \geq 1$ and $n \geq 2\ell + 4$.
\end{cor}

\noindent {\emph Proof}: Consider the poset $P_{\ell,n}$ whose elements are $p_i$ for $i \in \{1,\ldots,n\}$, where $p_i < p_j$ if and only if $i  + \ell < j$. We have $Av_{P_{\ell,n}}(132) \leq Av_{P_{\ell,n}}(123)$ with strict inequality if and only if $\ell \geq 1$ and $n \geq 2\ell + 4$ thanks to Theorem~\ref{thm:123}. Define the bijection~$\Phi\colon S_{P_{\ell,n}} \to S_n$ by  $\Phi(\sigma)_i = j$ if and only if~$\sigma_j = p_i$. Then a $132$ (respectively, $123$) pattern in $\sigma$ corresponds precisely to a $1-_\ell3-_\ell2$ (resp., $1-_\ell2-_\ell3$) gap pattern in $\Phi(\sigma)$. $\square$

These gap patterns apparently have not been introduced before but are similar in spirit to consecutive patterns. The $132$ versus $123$ dichotomy is present in consecutive pattern avoidance as well (see~\cite[Proposition 4.2]{elizalde}). However, it does not appear to us to be possible to recapture results about consecutive patterns using poset permutations.

\section{Open problems}

We hope that the simple, combinatorial proof of our main result will encourage further research into pattern avoidance in permutations on posets. Besides the chain and antichain (for which pattern avoidance is trivial), there are three other three-element posets underlying length three poset patterns: the ``hill" and ``valley" posets (two elements incomparable to each other, one element either greater than or less than both), and the ``line plus point" poset (two elements comparable, one incomparable to both). In the case of the hill and valley posets, we would like to see a proof of the following:

\begin{conj} \label{conj:main}
We have $Av_P(\{1\}\{1,2\}\{2\}) \leq Av_P(\{1\}\{2\}\{1,2\})$ for any poset $P$.\end{conj}

\noindent \emph{Remark}. Ideally one would find an injection from the set of $\{1\}\{1,2\}\{2\}$-avoiding permutations in $S_P$ to the set of $\{1\}\{2\}\{1,2\}$-avoiding permutations, perhaps using some subsequence similar to LRME. We have not been able to accomplish this. Both Simion and Schmidt's proof and the proof of our own Theorem~\ref{thm:123} rely on the fact that in any permutation of a totally ordered set (or poset), if there exists any 132 pattern then there exists some~132 pattern which begins with a left-to-right minimum (or minimal element). This is not true for $\{1\}\{2\}\{1,2\}$ and $\{1\}\{1,2\}\{2\}$. Moreover, it is not always possible to fix the LRME of a $\{1\}\{1,2\}\{2\}$-avoiding permutation and rearrange the remaining elements in a fashion which avoids $\{1\}\{2\}\{1,2\}$ as in the case of the permutation $(\{1\},\{2\},\{1,2\},\emptyset) \in S_{B_2}$. We might instead attempt to fix \emph{left-to-right minima} (LRM):  entries less than each preceding entry. But it is not always the case that if there exists a $\{1\}\{2\}\{12\}$ pattern then there exists such a pattern beginning with an LRM, as is evidenced by the permutation $(\{1,2,3\},\{1,2\},\{1,3\},\{1\},\{2\},\{3\},\{2,3\},\emptyset) \in S_{B_3}$. Nevertheless, computer tests on all posets with seven or fewer elements do suggest that for any set of LRM, there are at least as many ways to fill in the remaining elements, while preserving that set of LRM, and to avoid~$\{1\}\{2\}\{1,2\}$ as to avoid $\{1\}\{1,2\}\{2\}$. Therefore we suspect left-to-right minima may be a fruitful line of inquiry.

Furthermore, we expect strict inequality in Conjecture~\ref{conj:main} if and only if $P$ contains either of $R_1$ or $R_2$ below as an induced subposet:
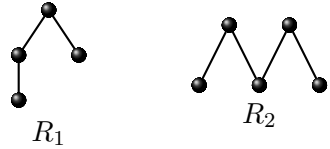
\begin{center}
\begin{tikzpicture}[scale=0.8]
\SetVertexMath
\GraphInit[vstyle=Art]
\SetUpVertex[MinSize=3pt]
\SetVertexLabel
\tikzset{VertexStyle/.style = {shape = circle,shading = ball,ball color = black,inner sep = 2pt}}
\SetUpEdge[color=black]

\Vertex[NoLabel,x=0,y=0.25]{a}
\Vertex[NoLabel,x=0,y=1]{b}
\Vertex[NoLabel,x=0.5,y=1.75]{c}
\Vertex[NoLabel,x=1,y=1]{d}
\Edges(a,b,c,d)

\draw (0.5,-0.3) node{$R_1$};

\Vertex[NoLabel,x=3,y=0.5]{a}
\Vertex[NoLabel,x=3.5,y=1.5]{b}
\Vertex[NoLabel,x=4,y=0.5]{c}
\Vertex[NoLabel,x=4.5,y=1.5]{d}
\Vertex[NoLabel,x=5,y=0.5]{e}
\Edges(a,b,c,d,e)

\draw (4,0) node{$R_2$};
\end{tikzpicture}
\end{center}
Computer tests on all posets $P$ with seven or fewer elements indicate that containment of either of $R_1$ and $R_2$ as an induced subposet is equivalent to the inequality between $Av_P(\{1\}\{1,2\}\{2\})$ and $Av_P(\{1\}\{2\}\{1,2\})$ being strict.  But even the ``if'' direction of this claim is more difficult than in the $Av_P(132) \leq Av_P(123)$ case, since it is not obvious, for example, which elements in $R_1$ or $R_2$ must be LRM. $\square$

The relationships between avoidance of $\{1\}\{1,2\}\{3\}$ and $\{1\}\{3\}\{1,2\}$, patterns arising from the final non-trivial three element poset, is more complicated. With posets $T$ and $U$ as below,
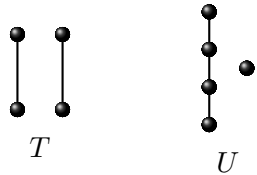
\begin{center}
\begin{tikzpicture}
\SetVertexMath
\GraphInit[vstyle=Art]
\SetUpVertex[MinSize=3pt]
\SetVertexLabel
\tikzset{VertexStyle/.style = {shape = circle,shading = ball,ball color = black,inner sep = 2pt}}
\SetUpEdge[color=black]

\Vertex[NoLabel,x=0.2,y=0.2]{a}
\Vertex[NoLabel,x=0.2,y=1.2]{b}
\Vertex[NoLabel,x=0.8,y=0.2]{c}
\Vertex[NoLabel,x=0.8,y=1.2]{d}
\Edges(a,b)
\Edges(c,d)

\draw (0.5,-0.3) node{$T$};

\Vertex[NoLabel,x=2.75,y=0]{a}
\Vertex[NoLabel,x=2.75,y=0.5]{b}
\Vertex[NoLabel,x=2.75,y=1]{c}
\Vertex[NoLabel,x=2.75,y=1.5]{d}
\Vertex[NoLabel,x=3.25,y=0.75]{e}
\Edges(a,b,c,d)

\draw (3,-0.5) node{$U$};
\end{tikzpicture}
\end{center}
we have 
\[Av_T(\{1\}\{1,2\}\{3\}) < Av_T(\{1\}\{3\}\{1,2\})\] 
but 
\[Av_U(\{1\}\{3\}\{1,2\}) < Av_U(\{1\}\{1,2\}\{3\}).\]
Considering these conflicting inequalities, we currently have no conjecture involving these patterns.

%T 14 < 12,  U 64 < 65

\noindent {\bf Acknowledgments}: This research was done at the East Tennessee State University Research Experience for Undergraduates with the support of NSF grant 1004624. We especially thank Professor Anant Godbole, director of the ETSU REU, for suggesting the area of investigation. We also thank Virginia Hogan, David Perkinson, and Yevgeniy Rudoy for their helpful comments.  All computer tests in our research were carried out using Sage mathematical software~\cite{sage-combinat}. Finally, we thank the anonymous referee for many helpful comments.

\end{document}